\documentclass[11pt,a4paper]{article}
\pdfoutput=1
\title{Natural $G$-Constellation Families}
\author{Timothy Logvinenko}
\usepackage{amsmath,amsfonts,amssymb,amsthm,epsfig,amscd,psfrag,latexsym,
comment}
\usepackage{graphicx}
\usepackage[all]{xy}
\usepackage[hyperref, colorlinks=true, pdfpagemode=none,
pdfmenubar=false, pdfstartview=FitH, linkcolor=blue, citecolor=blue,
urlcolor=blue, pdffitwindow=false]{hyperref}
\DeclareMathOperator{\homm}{Hom}
\DeclareMathOperator{\eend}{End}
\DeclareMathOperator{\seend}{\it \mathcal{E}nd}

\DeclareMathOperator{\gl}{GL}
\DeclareMathOperator{\gsl}{SL}

\DeclareMathOperator{\picr}{Pic}

\DeclareMathOperator{\divr}{Div}
\DeclareMathOperator{\cl}{Cl}
\DeclareMathOperator{\car}{Car}

\DeclareMathOperator{\spec}{Spec\;}

\DeclareMathOperator{\hilb}{Hilb}

\DeclareMathOperator{\ext}{Ext}

\DeclareMathOperator{\supp}{Supp}
\DeclareMathOperator{\ann}{Ann}

\DeclareMathOperator{\val}{\it val}
\DeclareMathOperator{\fract}{frac}
\DeclareMathOperator{\ramm}{Ram}

\DeclareMathOperator{\coh}{\bf Coh}
\DeclareMathOperator{\modd}{\bf Mod}
\begin{document}
\def\bv{\mathbf{v}}
\def\kgc_{K^*_G(\mathbb{C}^n)}
\def\kgchi_{K^*_\chi(\mathbb{C}^n)}
\def\kgcf_{K_G(\mathbb{C}^n)}
\def\kgchif_{K_\chi(\mathbb{C}^n)}
\def\gpic_{G\text{-}\picr}
\def\gcl_{G\text{-}\cl}
\def\trch_{{\chi_{0}}}
\def\regring{{R}}
\def\regrep{{V_{\text{reg}}}}
\def\givrep{{V_{\text{giv}}}}
\def\lbar{{(\mathbb{Z}^n)^\vee}}
\def\genpx_{{p_X}}
\def\genpy_{{p_Y}}
\def\genpcn_{p_{\mathbb{C}^n}}
\def\gnat{gnat}
\def\twalg{{\regring \rtimes G}}
\theoremstyle{definition}
\newtheorem{defn}{Definition}[section]
\newtheorem*{defn*}{Definition}
\newtheorem{exmpl}[defn]{Example}
\newtheorem*{exmpl*}{Example}
\newtheorem{exrc}[defn]{Exercise}
\newtheorem*{exrc*}{Exercise}
\newtheorem*{chk*}{Check}
\newtheorem*{remarks*}{Remarks}
\theoremstyle{plain}
\newtheorem{theorem}{Theorem}[section]
\newtheorem*{theorem*}{Theorem}
\newtheorem{conj}[defn]{Conjecture}
\newtheorem{prps}[defn]{Proposition}
\newtheorem*{prps*}{Proposition}
\newtheorem{cor}[defn]{Corollary}
\newtheorem*{cor*}{Corollary}
\newtheorem{lemma}[defn]{Lemma}
\newtheorem*{claim*}{Claim}
\numberwithin{equation}{section}
\maketitle

\begin{abstract}
Let $G$ be a finite subgroup of $\gl_n(\mathbb{C})$.
$G$-constellations are a scheme-theoretic generalization of orbits 
of $G$ in $\mathbb{C}^n$. We study flat families of $G$-constellations 
parametrised by an arbitrary resolution of the quotient space $\mathbb{C}^n/G$.
We develop a geometrical naturality criterion for such families, and
show that, for an abelian $G$, the number of equivalence classes of
these natural families is finite.

The main intended application is the derived McKay correspondence.

\end{abstract}

\setcounter{section}{-1} 

\section{Introduction} \label{section-intro}

Let $G \subseteq \gl_n(\mathbb{C})$ be a finite subgroup. In this
paper, we classify flat families of \it $G$-constellations \rm parametrised by 
a given resolution $Y$ of the singular quotient space $X =
\mathbb{C}^n / G$.
\begin{align*}
\xymatrix{ 
Y \ar[dr]^{\pi} & & \mathbb{C}^n \ar[dl]_{q} \\
& X &
}
\end{align*}

A $G$-constellation is a scheme-theoretical generalization of
a set-theoretical orbit of $G$ in $\mathbb{C}^n$. They first arose in 
the context of moduli space constructions of crepant resolutions 
of $X$. Interpreting $G$-constellations in terms of 
representations of the McKay quiver of $G$, it is possible to use
the methods of \cite{King94} to consruct via GIT fine moduli spaces 
of stable $G$-constellations. The main irreducible component of
such a moduli space turns out to be a projective crepant resolution of
$X$. By varying the stability parameter $\theta$ it is possible to
obtain different resolutions $M_\theta$. In case of $n = 3$ and $G$
abelian, it is possible to obtain all projective crepant resolutions in
this way \cite{Craw-Ishii-02}.  For further details see
\cite{Craw-thesis}, \cite{Craw-Ishii-02},
\cite{Craw-Maclagan-Thomas-05-I}, \cite{Craw-Maclagan-Thomas-05-II}.

The formal definition  of a $G$-constellation is:
 
\begin{defn}[\cite{Craw-thesis}]
A $G$-constellation is a $G$-equivariant coherent sheaf $\mathcal{F}$ 
whose global sections $\Gamma(\mathbb{C}^n, \mathcal{F})$, as a
representation of $G$, are isomorphic to the regular representation.  
\end{defn}

Families of $G$-constellations also occur naturally as objects defining 
\it Fourier-Mukai transforms \rm (cf.
\cite{BKR01},\cite{Craw-Ishii-02}, \cite{BonOr95}
\cite{Bridg97}) which give a category equivalence 
$D(Y) \xrightarrow{\sim} D^G(\mathbb{C}^n)$ between the bounded derived 
categories of coherent sheaves on $Y$ and of $G$-equivariant coherent 
sheaves on $\mathbb{C}^n$, respectively. This equivalence is 
known as the \it derived McKay correspondence \rm (cf.
\cite{Kinosaki-97}, \cite{BKR01},
\cite{Kawamata-LogCrepantBirationalMapsAndDerivedCategories},
\cite{Kaledin-05}). It is the derived 
category interpretation of the classical McKay correspondence between 
the representation theory of $G$ and the geometry of crepant
resolutions of $\mathbb{C}^n/G$. It was conjectured by Reid in 
\cite{Kinosaki-97} to hold for any finite subgroup $G$ of 
$\gsl_n(\mathbb{C})$ and any crepant resolution $Y$ of $\mathbb{C}^n/G$. 

In this paper we take an arbitrary resolution $Y \rightarrow
\mathbb{C}^n/G$ and prove that
it can support only a finite number (up to a twist by a line bundle) of
flat families of $G$-constellations.  We give a complete
classification of these families which allows one to explicitly
compute them. For the precise statement of the classification see the
end of this introduction.  

A motivation for this study is the fact that if a flat family of $G$-constellations on a crepant
resolution $Y$ of $\mathbb{C}^n/G$ is sufficiently orthogonal, 
then it defines an equivalence $D(Y) \rightarrow D^G(\mathbb{C}^n)$  
(\cite{Logvinenko-DerivedMcKayCorrespondenceViaPureSheafTransforms}, 
Theorem 1.1), i.e. the derived McKay correspondence conjecture
holds for $Y$. For an example of a specific application of this
see \cite{Logvinenko-DerivedMcKayCorrespondenceViaPureSheafTransforms}, 
\S 4, where the first known example of a derived McKay correspondence
for a non-projective crepant resolution is explicitly constructed.  

  This paper is laid out as follows. At the outset we allow $Y$ to 
be any normal scheme birational to the quotient space $X$ and 
first of all we move from the category
$\coh^G(\mathbb{C}^n)$ to the equivalent category 
$\modd^{\text{fg}}$-$\twalg$ of the finitely-generated modules
for the cross product algebra $\twalg$, where $R$ denotes the coordinate 
ring $\mathbb{C}[x_1, \dots, x_n]$ of $\mathbb{C}^n$. This makes
a family of $G$-constellations into a vector bundle on $Y$.
In Section \ref{section-gnatfam} we develop a geometrical naturality 
criterion for such families: mimicking the moduli spaces $M_\theta$ 
of $\theta$-stable $G$-constellations and their tautological
families, we demand for a $G$-constellation parametrised in a family
$\mathcal{F}$ by a point $p \in Y$ to be supported precisely on the
$G$-orbit corresponding to the point $\pi(p)$ in the quotient space
$X$. In other words, the support of the corresponding sheaf on $Y \times
\mathbb{C}^n$ must lie within the fibre product $Y \times_X
\mathbb{C}^n$. We call the families which satisfy this condition
$\gnat$-families (short for a \em geometrically natural\rm) and
demonstrate (Proposition \ref{prps-tfae}) that they enjoy a number of
other natural properties, including being equivalent (locally
isomorphic) to the natural family $\pi^* q_* \mathcal{O}_{\mathbb{C}^n}$
on the open set of $Y$ which lies over the free orbits in $X$. In this 
natural family a $G$-constellation which lies over a free orbit is 
the unique $G$-constellation supported on that orbit - its reduced subscheme
structure. Thus, in a sense, $\gnat$-families can be viewed as flat
deformations of free orbits of $G$. 
 
 Another property which characterises $\gnat$-families is that
it is possible to embed them into $K(\mathbb{C}^n)$, considered
as a constant sheaf on $Y$. This leads us to study $G$-equivariant 
locally free sub-$\mathcal{O}_Y$-modules of $K(\mathbb{C}^n)$. In Section 
\ref{section-valuations}, we study the rank one case. A $G$-invariant
invertible sub-$\mathcal{O}_Y$-module of $K(\mathbb{C}^n)$ is just 
a Cartier divisor, and we define $G\text{-}\car(Y)$, a group of 
$G$-Cartier divisors on $Y$, as a natural extension of the group 
of Cartier divisors which fits into a short exact sequence
\begin{align*}
 1 \rightarrow \car(Y) \rightarrow G\text{-}\car(Y) \xrightarrow{\rho} G^\vee \rightarrow 1
\end{align*}
where $G^\vee$ is the group of $1$-dimensional irreducible representations 
of $G$. 

 We then define $\mathbb{Q}$-valued valuations of these $G$-Cartier
divisors at prime Weil divisors of $Y$ and define $G$-$\divr Y$,
the group of $G$-Weil divisors of $Y$, as a torsion-free subgroup of 
$\mathbb{Q}$-Weil divisors which fits into a following exact sequence:
\begin{align*}
\xymatrix{ 
1 \ar[r] & \car Y \ar@{^{(}->}[d]_{\val_K} \ar[r] & G\text{-}\car Y \ar[d]_{\val_{K_G}}
\ar[r]^{\rho} & G^\vee \ar@{->>}[d]_{\val_{G^\vee}} \ar[r] & 1 \\
0 \ar[r] & \divr Y \ar[r] & G\text{-}\divr Y \ar[r] & 
\val_{G^\vee} (G^\vee)  \ar[r] & 0
}
\end{align*}
We then show that the three vertical maps in this diagram, $\val_K$,
the ordinary $\mathbb{Z}$-valued valuation of Cartier divisors, 
$\val_{K_G}$, the $\mathbb{Q}$-valued valuation of $G$-Cartier
divisors, and their quotient $\val_{G^\vee}$, a 
$\mathbb{Q}/\mathbb{Z}$-valued valuation of $G^\vee$, are all
isomorphisms when $Y$ is smooth and proper over $X$. 

Then, in Section \ref{section-classification}, we observe that when 
our group $G$ is abelian all its irreducible representations are of rank 
$1$, so any $\gnat$-family splits into invertible $G$-eigensheaves.
Thus $G$-Weil divisors are all that we need to classify it 
after an embedding into $K(\mathbb{C}^n)$. We further show that, since any 
$\gnat$-family $\mathcal{F}$ embedded into $K(\mathbb{C}^n)$ must 
be closed under the natural action of $\regring$ on the latter, all 
the $G$-eigensheaves into which $\mathcal{F}$ decomposes must be, in 
a certain sense, close to each other inside $K(\mathbb{C}^n)$. Up to a
twist by a line bundle, this leaves only a finite number of
possibilities for the corresponding $G$-Weil divisors.  Thus,
surprisingly, the number of equivalence classes of $\gnat$-families on
any $Y$ is finite. 

Our main result (Theorem \ref{theorem-classification}) is:
\begin{theorem*}[Classification of $\gnat$-families]

Let $G$ be a finite abelian subgroup of $\gl_n(\mathbb{C})$, $X$ the
quotient of $\mathbb{C}^n$ by the action of $G$ and $Y$ a 
resolution of $X$. Then isomorphism classes of \gnat-families on $Y$
are in $1$-to-$1$ correspondence with linear equivalence classes 
of $G$-divisor sets $\{D_\chi\}_{\chi \in G^\vee}$, each 
$D_\chi$ a $\chi$-Weil divisor, which satisfy the inequalities
\begin{align*}
D_\chi + (f) - D_{\chi \rho(f)} \geq 0 \quad \forall\; \chi \in
G^\vee, \text{$G$-homogeneous } f \in \regring
\end{align*}
Here $\rho(f) \in G^\vee$ is the homogeneous weight of $f$. 
Such a divisor set $\{D_\chi\}$ corresponds then to a 
$\gnat$-family $\bigoplus \mathcal{L}(-D_\chi)$.

This correspondence descends to a $1$-to-$1$ correspondence between 
equivalence classes of \gnat-families and sets $\{D_\chi\}$ 
as above and with $D_{\trch_} = 0$, where $\trch_$ is the trivial
character. Furthermore, each divisor $D_\chi$ in such a set 
satisfies
\begin{align*} M_\chi \geq
D_\chi \geq
- M_{\chi^{-1}}
\end{align*}
where $\{ M_{\chi} \}$ is a fixed divisor set defined by
\begin{align*}
M_\chi = \sum_P (\min_{f \in \regring_\chi} v_P(f)) P
\end{align*}
As a consequence, the number of equivalence classes of $\gnat$-families 
on $Y$ is finite. 
\end{theorem*}

\bf Acknowledgements: \rm The author would like to express his gratitude 
to Alastair Craw, Akira Ishii and Dmitry Kaledin for useful discussions 
on the subject and to Alastair King for the motivation, the discussions and 
the support. This paper was completed during the author's stay at RIMS, Kyoto,
and one would like to thank everyone at the institute for their hospitality. 

\section{$\gnat$-Families} \label{section-gnatfam}

\subsection{Families of $G$-Constellations}

Let $G$ be a finite abelian group and let $\givrep$ be an
$n$-dimensional faithful representation of $G$. We identify the symmetric 
algebra $S(\givrep^\vee)$ with the coordinate ring $\regring$ 
of $\mathbb{C}^n$ via a choice of such an isomorphism that 
the induced action of $G$ on $\mathbb{C}^n$ is diagonal. 
The (left) action of $G$ on $\givrep$ induces a (left) action of 
$G$ on $\regring$, where we adopt the convention that
\begin{align} \label{eq-gaction}
  g.f(\bv) = f(g^{-1}.\bv)\quad\quad\forall\; g \in G, f \in \regring,
\bv \in \givrep, \end{align}
When we consider the induced scheme morphisms
$g: \mathbb{C}^n \rightarrow \mathbb{C}^n$ and the induced sheaf
morphisms $g: \mathcal{O}_{\mathbb{C}^n} \rightarrow g^{-1}_{*}
\mathcal{O}_{\mathbb{C}_n}$, the convention above ensures that for any 
point $x \in \mathbb{C}^n$ and any function $f$ in the stalk
$\mathcal{O}_{\mathbb{C}^n,x}$ at $x$, the function $g.f$ is,
naturally, an element of the stalk $\mathcal{O}_{\mathbb{C}^n,g.x}$ at $g.x$

 Corresponding to the inclusion $\regring^G \subset \regring$ of the 
subring of $G$-invariant functions we have the quotient map $q:
\mathbb{C}^n \rightarrow X$, where $X = \spec \regring^G$ is 
the quotient space. This space is generally singular.

 We first wish to establish a notion of a family of \it $G$-Constellations \rm
parametrised by an arbitrary scheme.

\begin{defn}[\cite{Craw-Ishii-02}] \label{gcon-as-sheaf}
A \tt $G$-constellation \rm is a $G$-equivariant coherent sheaf
$\mathcal{F}$ on
$\mathbb{C}^n$ such that $H^0(\mathcal{F})$ is isomorphic, as a
$\mathbb{C}[G]$-module, to the regular representation $\regrep$. 
\end{defn}

 We would like for a family of $G$-constellations to be a locally
free sheaf on $Y$, whose restriction to any point of $Y$ would give
us the respective $G$-constellation. We'd like this restriction to 
be a finite-dimensional vector-space, and for this purpose, 
it would be better to consider, instead of the whole $G$-constellation
$\mathcal{F}$, just its space of global sections $\Gamma(\mathcal{F})$.
It is a vector space with $G$ and $\regring$ actions, satisfying
\begin{align} \label{o-g-equiv}
g.(f.\bv) = (g.f).(g.\bv)
\end{align}

On the other hand, for any vector space $V$ with $G$ and $\regring$
actions satisfying $\eqref{o-g-equiv}$, we can define maps $g:
\tilde{V} \rightarrow g^{-1}_* \tilde{V}$ to give the sheaf
$\tilde{V} = V \otimes_\regring \mathcal{O}_{\mathbb{C}^n}$ 
a $G$-equivariant structure. It is convenient to view such 
vector spaces as modules for the following non-commutative algebra:

\begin{defn} \label{defn-r-smash-g}
A cross-product algebra $\twalg$ is an algebra, which has the vector 
space structure of $\regring \otimes_{\mathbb{C}} \mathbb{C}[G]$ and 
the product defined by setting, for all $g_1, g_2 \in G$ and $f_1,
f_2 \in \regring$, 
\begin{align}
(f_1 \otimes g_1) \times (f_2 \otimes g_2) = (f_1 (g_1.f_2))
\otimes (g_1 g_2)
\end{align}
\end{defn}

This is not a pure formalism - $\twalg$ is one of the 
\em non-commutative crepant resolutions \rm of $\mathbb{C}^n / G$, 
a certain class of non-commutative algebras introduced by Michel van 
den Bergh in \cite{vdBergh2002} as an analogue of a commutative crepant
resolution for an arbitrary non-quotient Gorenstein singularity. 
For three-dimensional terminal singularities, van den Bergh shows
(\cite{vdBergh2002}, Theorem 6.3.1) that if a non-commutative 
crepant resolution $Q$ exists, then it is possible to construct
commutative crepant resolutions as moduli spaces of certain stable 
$Q$-modules.  

Functors $\Gamma(\bullet)$ and $\tilde{\bullet} = (\bullet)
\otimes_\regring \mathcal{O}_{\mathbb{C}^n}$ give an equivalence
(compare to \cite{Harts77}, p. 113, Corollary 5.5) between the 
categories of quasi-coherent $G$-equivariant sheaves on $\mathbb{C}^n$ 
and of $\twalg$-modules. $G$-constellations then 
correspond to $\twalg$-modules, whose underlying
$G$-representation is $\regrep$. As an abuse of notation, we shall 
use the term `$G$-constellation' to refer to both the equivariant
sheaf and the corresponding $\twalg$-module. 

\begin{defn} \label{def-gcon-fam}
A \tt family of $G$-constellations parametrised by a scheme 
$S$ \rm is a sheaf $\mathcal{F}$ of 
$(\twalg) \otimes_\mathbb{C} \mathcal{O}_S$-modules, locally
free as an $\mathcal{O}_S$-module, and such that for any 
point $\iota_p:\; \spec \mathbb{C} \hookrightarrow S$,
its fiber $\mathcal{F}_{|p} = \iota_p^* \mathcal{F}$ 
is a $G$-constellation.

We shall say that two families $\mathcal{F}$ and $\mathcal{F}'$ are 
equivalent if they are locally isomorphic as $(\twalg)
\otimes_\mathbb{C} \mathcal{O}_S$-modules.
\end{defn}

\subsection{$\gnat$-Families}

  Let $Y$ be a normal scheme and $\pi: Y \rightarrow X$ be a 
birational map. 

\begin{align*}
\xymatrix{ 
Y \ar[dr]^{\pi} & & \mathbb{C}^n \ar[dl]_{q} \\
& X &
}
\end{align*}

  We wish to refine the definition
$\eqref{def-gcon-fam}$ above and develop a notion of a geometrically
natural family of $G$-constellations parametrised by $Y$.  

  Any free $G$-orbit supports a unique $G$-cluster $Z \subset \mathbb{C}^n$:
the reduced induced closed subscheme structure. Let $U$ be an open subset 
of $Y$ such that $\pi(U)$ consists of free orbits of $G$ and consider the
sheaf $\pi^* q_* \mathcal{O}_{\mathbb{C}^n}$ restricted to $U$. It
has a natural $(\twalg)$-module structure 
induced from $\mathcal{O}_{\mathbb{C}^n}$. It is locally free as an 
$\mathcal{O}_U$ module, since the quotient map $q$ is flat wherever
$G$ acts freely. Its fiber at a point $p \in Y$ is $\Gamma(\mathcal{O}_Z)$, 
where $Z$ is the $G$-cluster corresponding to the free orbit $q^{-1}\pi(p)$. 
Thus $\pi^* q_* \mathcal{O}_{\mathbb{C}^n}$ is a natural family 
of $G$-constellations, indeed of $G$-clusters, on $U \subset Y$. 

Its fiber at the generic point of $Y$ is $K(\mathbb{C}^n)$. The 
Normal Basis Theorem from Galois theory (\cite{Garl86}, Theorem 19.6)
gives an isomorphism from $K(\mathbb{C}^n)$ to the generic fiber of 
any $G$-constellation family on $Y$, which we can write as $K(Y)
\otimes_\mathbb{C} \regrep$, but this isomorphism is only 
$K(Y)$ and $G$, but not necessarily $\regring$, equivariant.  

On the other hand, for any $G$-constellation in a sense of 
$G$-equivariant sheaf, we can consider its support in $\mathbb{C}^n$. 
For instance, in the natural family $\pi^* q_* \mathcal{O}_{\mathbb{C}^n}$
discussed above the support of the $G$-constellation parametrised by a 
point $p \in U$ is precisely the $G$-orbit $q^{-1}\pi(p)$. 
This turns out to be the criterion we seek and we shall
show that any family satisfying it is generically equivalent to the 
natural one.

\begin{defn}\label{defn-gnat-family}
A \tt $\gnat$-family $\mathcal{F}$ \rm (short for \em geometrically 
natural \rm family) is a family of $G$-constellations parametrised 
by $Y$ such that for any $p \in Y$
\begin{align}\label{eqn-support-equality}
q \left(\supp_{\mathbb{C}^n} \mathcal{F}_{|p}\right) = \pi(p)
\end{align}
\end{defn}

\begin{prps} \label{prps-tfae}
Let $Y$ be a normal scheme and $\pi: Y \rightarrow X$ be a 
birational map. 
Let $\mathcal{F}$ be a family of $G$-constellations on $Y$.  
Then the following are equivalent:
\begin{enumerate}
\item \label{item-local-equival} On any $U \subset Y$, such that $\pi U$ consists of 
free orbits, $\mathcal{F}$ is equivalent to 
$\pi^* q_* \mathcal{O}_{\mathbb{C}^n}$.
\item \label{item-generically-natural} There exists an $(\twalg) \otimes_\mathbb{C} K(Y)$-module isomorphism:
\begin{align*}
\mathcal{F}_{|p_Y} \xrightarrow{\sim} 
(\pi^* q_* \mathcal{O}_{\mathbb{C}^n})_{p_Y} 
\end{align*}
where $p_Y$ is the generic point of $Y$. 
\item \label{item-embeds-into-kcn} There exists an $(\twalg) \otimes_\mathbb{C} \mathcal{O}_Y$-module
embedding $$F \hookrightarrow K(\mathbb{C}^n)$$ where
$K(\mathbb{C}^n)$ is viewed as a constant sheaf on $Y$ and
given a $\mathcal{O}_Y$-module structure via the birational map 
$\pi: Y \rightarrow X$.
\item \label{item-gnat-family} $\mathcal{F}$ is a \gnat-family. 
\item \label{item-action-descent} The action of $(\twalg)
\otimes_\mathbb{C} \mathcal{O}_Y$ on $\mathcal{F}$ descends to the
action of $(\twalg)~\otimes_{\regring^G }~\mathcal{O}_Y$, where 
$\regring^G$-module structure on $\mathcal{O}_Y$ is induced by the map
$\pi: Y \rightarrow X$.
\end{enumerate}
\end{prps}

\begin{proof}
$\ref{item-local-equival} \Rightarrow \ref{item-generically-natural}$
is restricting any of the local isomorphisms to the stalk at the
generic point $p_Y$ of $Y$. $\ref{item-generically-natural} \Rightarrow
\ref{item-embeds-into-kcn}$: the embedding is given by 
the natural map $\mathcal{F} \hookrightarrow \mathcal{F} \otimes K(Y)$.
As $Y$ is irreducible and $\mathcal{F}$ is locally free, 
$\mathcal{F}~\otimes~K(Y)$ is isomorphic to $\mathcal{F}_{p_Y}$, and
hence to $K(\mathbb{C}^n)$. $\ref{item-embeds-into-kcn} \Rightarrow
\ref{item-action-descent}$ is immediate by inspecting the natural 
$\twalg \otimes_{\mathbb{C}} \mathcal{O}_Y$-module structure on
$K(\mathbb{C}^n)$. $\ref{item-action-descent} \Rightarrow 
\ref{item-gnat-family}$ is also immediate, as the descent of the
action of $\twalg \otimes_\mathbb{C} \mathcal{O}_Y$ to that
of $\twalg \otimes_{\regring^G} \mathcal{O}_Y$ implies that for 
any $p \in Y$ we have 
$\mathfrak{m}_{\pi(p)} \subset \ann_{\regring} \mathcal{F}_{|p}$, 
where $\mathfrak{m}_{\pi(p)} \subset \regring^G$ is the maximal ideal 
of $\pi(p)$. Therefore $\mathfrak{m}_{\pi(p)} = (\ann_{\regring} \mathcal{F}_{|p})^G$, which is equivalent to 
$\eqref{eqn-support-equality}$.

$\ref{item-gnat-family} \Rightarrow \ref{item-action-descent}$:
Consider the following composition of algebra morphisms:
\begin{align*}
\twalg \otimes_{\mathbb{C}} \mathcal{O}_Y \xrightarrow{\alpha}
\seend_{\mathcal{O}_Y}(\mathcal{F}) \xrightarrow{\beta_p}
\eend_\mathbb{C}(\mathcal{F}_{|p})
\end{align*}
where $\alpha$ is the action map of $\twalg \otimes_{\mathbb{C}}
\mathcal{O}_Y$ on $\mathcal{F}$ and $\beta_p$ is restriction 
to the fiber at a point $p \in Y$.

To show that $\alpha$ filters through 
$\twalg \otimes_{\regring^G} \mathcal{O}_Y$ it suffices to show
that for any $f \in \regring^G$ we have $ f \otimes 1 - 1 \otimes f
\in \ker(\alpha)$. From $\eqref{eqn-support-equality}$ we have $\mathfrak{m}_{\pi(p)} = (\ann_{\regring}
\mathcal{F}_{|p})^G$, and therefore
$$ \beta_p \alpha ( (f - f(p)) \otimes 1 ) = 0 $$
Observe that $ \beta_p \alpha (f(p) \otimes 1) = 
f(p)\; 1_{\eend_\mathbb{C} \mathcal{F}_{|p}} = 
\beta_p \alpha (1 \otimes f)$, and therefore
\begin{align}\label{eqn-rg-action-commuting}
\beta_p \alpha ( f \otimes 1 - 1 \otimes f ) = 0
\end{align}
As $\seend_{\mathcal{O}_Y} \mathcal{F}$ is locally free, 
$\eqref{eqn-rg-action-commuting}$ holding $\forall\; p \in Y$ implies 
$\alpha(f \otimes 1 - 1 \otimes f) = 0$, as required.

$\ref{item-action-descent} \Rightarrow \ref{item-local-equival}$:
We have the $\twalg
\otimes_{\regring^G} \mathcal{O}_Y$-action on $\mathcal{F}$:
$$ 
\quad \twalg \otimes_{\regring^G} \mathcal{O}_Y 
\xrightarrow{\alpha}
\seend_{\mathcal{O}_Y}(\mathcal{F}) 
$$
LHS is isomorphic to $\pi^* \seend_{\mathcal{O}_X}(q_*
\mathcal{O}_{\mathbb{C}^n})$. Over $U$, since $q$ is 
flat over $\pi(U)$, LHS is further isomorphic to
$\seend_{\mathcal{O}_U}(\pi^*
q_* \mathcal{O}_{\mathbb{C}_n})$. Thus we have:
\begin{align}\label{eqn-endom-level-action}
\seend_{\mathcal{O}_U}(\pi^* q_* \mathcal{O}_{\mathbb{C}_n})
\xrightarrow{\alpha'}
\seend_{\mathcal{O}_U}(\mathcal{F}) 
\end{align}
This map $\eqref{eqn-endom-level-action}$ is an
$\mathcal{O}_U$-algebra homomorphism of (split) Azumaya algebras over
$U$ of the same rank. By a general result on Azumaya algebras 
any such is an isomorphism
(see \cite{Adjamagbo-Charbonnel-VanDenEssen-05}, Theorem 5.3, for full
generality, but the original result in \cite{Auslander-Goldamn-1960},
Corollary 3.4 will also suffice here). Now
Skolem-Noether theorem for Azumaya algebras
(\cite{Milne-EtaleCohomology}, IV, \S2, Proposition 2.3) implies
that locally $\alpha'$ must be induced by isomorphisms $\pi^* q_*
\mathcal{O}_{\mathbb{C}^n} \xrightarrow{\sim} \mathcal{F}$. 
\end{proof}

\section{$G$-Cartier and $G$-Weil divisors} \label{section-valuations}

If $\mathcal{F}$ is a $\gnat$-family, by Proposition $\ref{prps-tfae}$
we can embed it into $K(\mathbb{C}^n)$. We need, therefore, to study
$G$-subsheaves of $K(\mathbb{C}^n)$ which are locally free on $Y$. In
this section we treat the rank 1 case, i.e. the invertible sheaves. 
Now, on an arbitrary scheme $S$, an invertible sheaf
together with its embedding into $K(S)$ defines a unique Cartier
divisor on $S$. But here we embed not into $K(Y)$ but into its
Galois extension $K(\mathbb{C}^n)$. Recall that we identify $K(Y)$
with $K(\mathbb{C}^n)^G$ via the birational map $Y \xrightarrow{\pi}
X$.
We therefore seek to extend 
the familiar construction of Cartier divisors to accommodate for 
this fact.

\subsection{G-Cartier divisors} \label{subsection-gcartier}

We write $G^\vee$ for $\homm(G, \mathbb{C}^*)$, the group of 
irreducible representations of $G$ of rank 1.
\begin{defn}
 We shall say that a rational function $f \in K(\mathbb{C}^n)$ is \tt
$G$-homogeneous of weight $\chi \in G^\vee$ \rm if 
\begin{align} \label{eq-ghom-fn}
g.f = \chi(g^{-1}) f \quad \forall\; g \in G 
\end{align}

 We shall denote by $K_\chi(\mathbb{C}^n)$ the subset of $K(\mathbb{C}^n)$ of
homogeneous elements of a specific weight $\chi$ 
and by $K_G(\mathbb{C}^n)$ the subset of $K(\mathbb{C}^n)$ of
all the $G$-homogeneous elements. 
We shall use $\regring_\chi$ and $\regring_G$ to mean
$\regring \cap K_\chi(\mathbb{C}^n)$ and $\regring \cap
K_G(\mathbb{C}^n)$ respectively. 
\end{defn} 

\bf NB: \rm The choice of a sign is dictated by wanting $f \in \regring$ 
to be homogeneous of weight $\chi \in G^\vee$ if $f(g.v) = \chi(g) f(v)$ 
for all $g \in G$ and $v \in \mathbb{C}^n$. 

The invertible elements of
$K_G(\mathbb{C}^n)$ form a multiplicative group which we shall denote 
by $\kgc_$. We have a short exact sequence: 
\begin{align}\label{eqn-kgc-ses}
 1 \rightarrow K^*(Y) \rightarrow \kgc_ \xrightarrow{\rho} G^\vee \rightarrow 1
\end{align}

 The following replicates, almost word-for-word, the definition of
a Cartier divisor in \cite{Harts77}, pp. 140-141. 

\begin{defn}
  A group of \tt $G$-Cartier divisors \rm on $Y$, denoted by
$G$-$\car(Y)$ 
is the group of global sections of the sheaf of multiplicative groups 
$\kgc_ / \mathcal{O}^*_Y$, i.e. the quotient of the constant sheaf 
$\kgc_$ on $Y$ by the sheaf $\mathcal{O}^*_Y$ of invertible regular 
functions. 
\end{defn}

  Observe that $\eqref{eqn-kgc-ses}$ gives a well-defined short 
exact sequence:
\begin{align}\label{eqn-car-ses}
 1 \rightarrow \car(Y) \rightarrow G\text{-}\car(Y) \xrightarrow{\rho} G^\vee \rightarrow 1
\end{align}

Given a $G$-Cartier divisor, we call its image $\chi \in G^\vee$ under $\rho$ 
\tt the weight \rm of the divisor and say, further, that the divisor 
is \tt $\chi$-Cartier\rm. 

  A $G$-Cartier divisor can be specified by a choice of an open cover 
$\{U_i\}$ of $Y$ and functions $\{f_i\} \subseteq \kgc_$ such 
that $f_i / f_j \in \Gamma(U_i \cap U_j, \mathcal{O}^*_Y)$. 
In such case, the weight of the divisor is the weight of any one 
of $f_i$.   

 As with ordinary Cartier divisors, we say that a $G$-Cartier divisor is 
principal if it lies in the image of the natural map 
$\kgc_ \rightarrow \kgc_ / \mathcal{O}^*_Y$ and call two divisors 
linearly equivalent if their difference is principal. 

 Consider now a $\chi$-Cartier divisor $D$ on $Y$ specified 
by a collection $\{(U_i, f_i)\}$ where $U_i$ form an open cover 
of $Y$ and $f_i \in \kgchi_$.  
We define an invertible sheaf $\mathcal{L}(D)$ on $Y$ as
the sub-$\mathcal{O}_Y$-module of $K(\mathbb{C}^n)$ generated by
$f_i^{-1}$ on $U_i$. Observe that $G$ acts on
$\mathcal{L}(D)$, the action being the restriction of the one on 
$K(\mathbb{C}^n)$, and that it acts on every section by 
the character $\chi$.  

\begin{prps} \label{prps-gcar-gpic}
  The map $D \rightarrow \mathcal{L}(D)$ gives an isomorphism between 
$G$-$\car Y$ and the group of invertible $G$-subsheaves of
$K(\mathbb{C}^n)$. Furthermore, it descends to an isomorphism
of the group $\gcl_$ of $G$-Cartier divisors up to linear
equivalence and the group $\gpic_$ of invertible $G$-sheaves on $Y$.
\end{prps}
\begin{proof}
A standard argument in \cite{Harts77}, Proposition 6.13, shows
everything claimed, apart from the fact we can embed any invertible 
$G$-sheaf $\mathcal{L}$, with $G$ acting by some $\chi \in G^\vee$, 
as a sub-$\mathcal{O}_Y$-module into $K(\mathbb{C}^n)$. 
  
  Given such $\mathcal{L}$, we consider the sheaf 
$\mathcal{L} \otimes_{\mathcal{O}_Y} K(Y)$. On every open set $U_i$ where 
$\mathcal{L}$ is trivial, it is $G$-equivariantly isomorphic to the constant 
sheaf $K_\chi(\mathbb{C}^n)$. On an irreducible scheme a sheaf
constant on an open cover is constant itself, so as $Y$ is irreducible we have 
$\mathcal{L} \otimes_{\mathcal{O}_Y} K(Y) \simeq K_\chi(\mathbb{C}^n)$ and a 
particular choice of this isomorphism gives the necessary embedding as
\begin{align*}
\mathcal{L} \rightarrow \mathcal{L} \otimes_{\mathcal{O}_Y} K(Y)
\simeq K_\chi(\mathbb{C}^n) \subset K(\mathbb{C}^n)
\end{align*}
\end{proof}

\subsection{Homogeneous valuations}
\label{subsection-homogeneous-valuations}

 We now aim to develop a matching notion of $G$-Weil divisors. Recall
that the homomorphism from ordinary Cartier to ordinary Weil divisors
is defined in terms of valuations of rational functions at prime 
Weil divisors of $Y$. 

 Valuations at prime divisors of $Y$ define a unique group
homomorphism $\val_K$ from $K^*(Y)$ to $\divr Y$, the group 
of Weil divisors. Looking at the short exact sequence
$\eqref{eqn-kgc-ses}$, we see that $\val_K$ must extend uniquely
to a homomorphism $\val_{K_G}$ from $\kgc_$ to $\mathbb{Q}\text{-}\divr Y$,
as $G^\vee$ is finite and $\mathbb{Q}$ is injective. We further
obtain a quotient homomorphism $\val_{G^\vee}$ from $G^\vee$ to
$\mathbb{Q}/\mathbb{Z}\text{-}\divr Y$.  

 Explicitly, we set:
\begin{defn} \label{defn-homogeneous-valuations}
 Let $P$ be a prime Weil divisor on $Y$. 

 For any $f \in \kgc_$, observe that $f^{|G|}$ is necessarily 
of trivial weight and hence lies in $K(Y)$. We 
define \tt valuation of $f$ at $P$ \rm to be 
\begin{align}
v_P(f) = \frac{1}{|G|}v_P(f^{|G|}) \in \mathbb{Q}
\end{align}
where $v_P(f^{|G|})$ is the ordinary valuation in 
the local ring of $P$. 

 For any $\chi \in G^\vee$, observe that for any $f, f'$ homogeneous
of weight $\chi$ their ratio $f/f'$ is of trivial character and 
therefore has integer valuation. We define \tt valuation of 
$\chi$ at $P$ \rm to be
\begin{align}
v_P(\chi) = \fract( v_P(f)) \in \mathbb{Q}/\mathbb{Z}
\end{align}
 where $f$ is any homogeneous function of weight $\chi$ and
$\fract(\text{-})$ denotes the fractional part. 
\end{defn}

It can be readily verified that $\val_{K_G} = \sum v_P(\text{-}) P$
and $\val_{G^\vee} = \sum v_P(\text{-})P$. Furthermore,  
the short exact sequence $\eqref{eqn-car-ses}$ 
becomes a commutative diagram:
\begin{align}
\xymatrix{ 
1 \ar[r] & \car Y \ar[d]_{\val_K} \ar[r] & G\text{-}\car Y \ar[d]_{\val_{K_G}}
\ar[r]^{\rho} & G^\vee \ar[d]_{\val_{G^\vee}} \ar[r] & 1 \\
0 \ar[r] & \divr Y \ar[r] & \mathbb{Q}\text{-}\divr Y \ar[r] & 
\mathbb{Q}/\mathbb{Z}{-}\divr Y \ar[r] & 0
}
\end{align}

\subsection{$G$-Weil divisors}
\label{subsection-gweil-divisors}

Aiming to have a short exact sequence similar to $\eqref{eqn-car-ses}$, 
we now define the group $G$-$\divr Y$ 
of $G$-Weil divisors
to be the subgroup of $\mathbb{Q}$-$\divr Y$, which consists of the 
pre-images of 
$\val_{G^\vee} (G^\vee) \subset \mathbb{Q}/\mathbb{Z}$-$\divr Y$. 

\begin{defn} \label{defn-gweil-divisor}
We say that a $\mathbb{Q}$-Weil divisor $\sum q_P P$ on $Y$ is 
\tt a $G$-Weil divisor \rm if there exists $\chi \in G^\vee$ 
such that 
\begin{align} \label{eqn-gweil-cond}
\fract(q_P) = v_P(\chi) \quad \quad \text{ for all prime Weil } P
\end{align}
\end{defn}

We call a $G$-Weil divisor principal if it is an image of
a single function $f \in \kgc_$ under $\val_{K^G}$, call 
two $G$-Weil divisors linearly equivalent if their difference
is principal and call a divisor $\sum q_i D_i$ effective if
all $q_i \geq 0$. 

We now have a following commutative diagram:
\begin{align} \label{diag-gcar-gweil}
\xymatrix{ 
1 \ar[r] & \car Y \ar@{^{(}->}[d]_{\val_K} \ar[r] & G\text{-}\car Y \ar[d]_{\val_{K_G}}
\ar[r]^{\rho} & G^\vee \ar@{->>}[d]_{\val_{G^\vee}} \ar[r] & 1 \\
0 \ar[r] & \divr Y \ar[r] & G\text{-}\divr Y \ar[r] & 
\val_{G^\vee} (G^\vee)  \ar[r] & 0
}
\end{align}

 A warning: for general $Y$, even a smooth one, $G$-Cartier and
$G$-Weil divisors may not be very well behaved. For an example let $Y$
be the smooth locus of $X$.  It can be shown, that while $\val_K$ is
an isomorphism, $\val_{K_G}$ is not even injective as $G$-$\car Y$ has
torsion. And $\val_{G^\vee}$ is the zero map, thus $G$-$\divr Y$ is just
$\divr Y$. 

\begin{prps}\label{prps-resolution-G-cartier-equals-G-Weil}
If $Y$ is smooth and proper over $X$, then $\val_K$, $\val_{K_G}$
and $\val_{G^\vee}$ in $\eqref{diag-gcar-gweil}$ are isomorphisms.
\end{prps}

\begin{proof}
  If $Y$ is smooth, or at least locally factorial, $\val_K$ is well-known 
to be an isomorphism (\cite{Harts77}, Proposition 6.11). It therefore
suffices to show that $\val_{G^\vee}$ is injective and hence an isomorphism. 
As diagram $\eqref{diag-gcar-gweil}$ commutes, $\val_{K_G}$ will
then also have to be an isomorphism.

  Fix $\chi \in G^\vee$. Let $Y_\chi$ denote the normalisation of 
$Y \times_X (\mathbb{C}^n / \ker \chi)$. It is a Galois covering of $Y$
whose Galois group is $\chi(G)$. By Zariski-Nagata's purity
of the branch locus theorem (\cite{Zariski-1958}, Proposition 2), 
the ramification locus of $Y_\chi \rightarrow Y$ is either empty 
or of pure codimension one. As $Y$ is smooth, $Y_\chi \rightarrow Y$
being finite and unramified would make it an \'etale cover. Which is
impossible, since a resolution of a quotient singularity is well 
known to be simply-connected (see, for instance, \cite{Verbitsky-2000}, 
Theorem 4.1).

  Thus, we can assume there exists a ramification divisor 
$P \subset Y_\chi$. Let $Q$ be its image in $Y$.
Let $\ramm(P)$ be the subgroup of $G$ which fixes
$P$ pointwise. Then 
$n_{\small \text{ram}} = \left| \ramm(P) / \ker \chi\right|$ 
is the ramification index of $P$. We can take ordinary 
integer valuations of $\kgchi_$ on prime divisors of $Y_\chi$ as 
$\kgchi_ \subset K(\mathbb{C}^n)^{\ker \chi}$. It is easy to see that 
for any $f \in \kgchi_$ 
\begin{align} \label{eqn-ramification-valuation}
v_{Q}(f) = \frac{1}{n_{\text{\small ram}}}v_{P}(f)
\end{align}
where LHS is a rational valuation in sense of Definition 
\ref{defn-homogeneous-valuations}.

If $v_{Q}(\chi) = 0$, then $v_{Q}(\kgchi_) \subset \mathbb{Z}$. Then
necessarily $v_{Q}(\kgchi_) = \mathbb{Z}$, as $\kgchi_$ is a coset 
of $K(Y)$ in $\kgc_$. In particular, there would exist 
$f_\chi \in \kgchi_$, such that $v_{Q}(f_\chi) = 0$, i.e. $f_\chi$ 
is a unit in $\mathcal{O}_{Y_\chi, P}$.  
Which is impossible: any $g \in \ramm(P)$ fixes $P$ pointwise, 
in particular $f - g.f \in \mathfrak{m}_{Y, P}$ for any 
$f \in \mathcal{O}_{Y, P}$. As $\ramm(P) / \ker \chi$ is non-trivial
we can choose $g$ such that $\chi(g) \neq 1$ and then
$f_\chi = \frac{1}{1 - \chi(g)} (f_\chi - g.f_\chi)$ must lie in 
$\mathfrak{m}_{Y, P}$. This finishes the proof. 

For abelian $G$, this all can be seen very explicitly by exploiting 
the toric structure
of the singularity: even though we do not assume the resolution $Y$ to 
be toric, it has been proven by Bouvier
(\cite{Bouvier-1998}, Theorem 1.1) and by Ishii and Koll\'ar 
(\cite{Kollar-Ishii-2003}, Corollary 3.17, in a more general context 
of Nash problem) that every essential divisor over $X$ (i.e. a divisor
which must appear on every resolution) is toric. The set of 
essential toric divisors is well understood - it can be 
identified with the Hilbert basis of the positive octant of 
the toric lattice of weights, and then with a subset of $\ext^1(G^\vee, \mathbb{Z}) = \homm(G^\vee,
\mathbb{Q}/\mathbb{Z})$. This correspondence sends each divisor
precisely to the valuation of $G^\vee$ at it, see
\cite{Logvinenko-thesis}, Section 4.3 for more detail. 
\end{proof}

We also show that, away from a finite
number of prime divisors on $Y$, all $G$-Weil divisors are ordinary Weil. 

\begin{prps}
Unless a prime divisor $P \subset Y$ is exceptional or its image
in $X$ is a branch divisor of $\mathbb{C}^n \rightarrow X$, 
the valuation $v_P: G^\vee \rightarrow \mathbb{Q}/\mathbb{Z}$
is the zero-map. 
\end{prps}
\begin{proof}
If $P$ is not exceptional, let $Q$ be its image in $X$. The valuations
at $P$ and $Q$ are the same, so it suffices to prove the statement
about $v_Q$. Let $P'$ be any divisor in $\mathbb{C}^n$ which lies above 
$Q$. As in Proposition \ref{prps-resolution-G-cartier-equals-G-Weil}, 
for any $f \in \kgc_$ we have 
$v_{Q}(f) = \frac{1}{n_{\text{\small ram}}}v_{P'}(f)$ where 
$n_{\text{\small ram}}$ is the ramification index of $P'$. 
Unless $Q$ is a branch divisor, $n_{\text{\small ram}} = 1$
and $v_{Q} = v_{P'}$. Which makes $v_{Q}$ integer-valued on
$\kgc_$ and makes the quotient homomorphism 
$G^\vee \rightarrow \mathbb{Q}/\mathbb{Z}$ the zero map. 
\end{proof}

\section{Classification of $\gnat$-families} \label{section-classification}

\subsection{Reductor Sets}

 From now on, in addition to assuming that $G$ is a finite group 
acting faithfully on $\givrep$, we also assume that $G$ is abelian.
We further assume that $Y$ is smooth and $\pi: Y \rightarrow X$ is 
proper. 

 Let $\mathcal{F}$ be a \gnat-family on $Y$. Write the decomposition 
of $\mathcal{F}$ into $G$-eigensheaves as $\bigoplus_{\chi \in G^\vee}
\mathcal{F}_\chi$. By Proposition \ref{prps-tfae} we can  
embed $\mathcal{F}$ into $K(\mathbb{C}^n)$
and, as was demonstrated in Proposition \ref{prps-gcar-gpic}, 
the image of each $\mathcal{F}_\chi$ defines a $\chi$-Cartier 
divisor. Hence $\mathcal{F} \simeq \bigoplus_\chi \mathcal{L}( - D_\chi)$
for some set $\{D_\chi\}_{\chi \in G^\vee}$ of $G$-Weil divisors.

\begin{defn}
Let $\{D_\chi\}_{\chi \in G^\vee}$ be a set of $G$-Weil divisors on $Y$. 
We call it a \tt reductor set \rm if each $D_\chi$ is 
a $\chi$-Weil divisor and $\oplus \mathcal{L}(-D_\chi)$ 
is a \gnat-family on $Y$. We call a reductor set \tt normalised \rm 
if $D_{\trch_} = 0$. We say that two reductor sets $\{D_\chi\}$ and
$\{D'_\chi\}$ are linearly equivalent if there exists $f \in K(Y)$
such that $D_\chi - D'_\chi = \divr f$ for all $\chi \in G^\vee$.
\end{defn}

\begin{lemma} \label{lemma-morphism-gnatfam}
Let $\{D_\chi\}$ and $\{D'_\chi\}$ be two reductor sets.
Any $(\twalg) \otimes \mathcal{O}_Y$-module morphism $\phi:\; 
\bigoplus \mathcal{L}(-D_\chi) \rightarrow \bigoplus \mathcal{L}(-D'_\chi)$ 
is necessarily a multiplication inside $K(\mathbb{C}^n)$ by some $f \in K(Y)$. 
\end{lemma}
\begin{proof}
Because of $G$-equivariance $\phi$ decomposes as 
$\bigoplus_{\chi \in G^\vee} \phi_\chi$ with $\phi_{\chi}$ a morphism
$\mathcal{L}(-D_\chi)~\rightarrow~\mathcal{L}(-D'_\chi)$.
Each $\phi_{\chi}$ is a morphism of invertible sub-$\mathcal{O}_Y$-modules 
of $K(\mathbb{C}^n)$ and so is necessarily a multiplication by some
$f_\chi \in K(Y)$: consider the induced map $\mathcal{O}_Y
\rightarrow \mathcal{L}(-D_\chi + D'_\chi)$ and take $f_{\chi}$
to be the image of $1$ under this map. 

It remains to show that all $f_\chi$ are equal. Fix any $\chi \in G^\vee$ 
and consider any $G$-homogeneous $m \in R$ of weight $\chi$.
Take any $s \in \mathcal{L}(- D_\trch_)$. Then $m s \in \mathcal{L}(- D_\chi)$ 
and by $\regring$-equivariance of $\phi$
\begin{align}
\phi_\chi(ms) = m \phi_{\trch_}(s) = f_{\trch_} m s
\end{align}
and hence $f_{\chi} = f_{\trch_}$ for all $\chi \in G^\vee$. 
\end{proof}

\begin{cor} \label{cor-iso-families-reductors-lineq}
Isomorphism classes of $\gnat$-families on $Y$ are in 1-to-1
correspondence with linear equivalence classes of reductor sets. 
\end{cor}
\begin{proof}
If in the proof of Lemma $\ref{lemma-morphism-gnatfam}$ each $\phi_\chi$
is an isomorphism, then $f$, by construction, globally generates
each $\mathcal{L}(- D'_\chi + D_\chi)$. Thus $D_\chi - D'_\chi =
\divr(f)$. 
\end{proof}

\begin{prps}\label{prps-equiv-families-divisor-diff}
Let $\{D_\chi\}$ and $\{D'_\chi\}$ be two reductor sets.
Then $\bigoplus \mathcal{L}(-D_\chi)$ and $\bigoplus
\mathcal{L}(-D'_\chi)$ are equivalent (locally isomorphic)
if and only if there exists a Weil divisor $N$ such that
$D_\chi - D'_\chi = N$ for all $\chi \in G^\vee$. 
\end{prps}

\begin{proof}
The `if' direction is immediate.

Conversely, if the families are equivalent, then by applying Lemma
$\ref{lemma-morphism-gnatfam}$ to each local isomorphism, we 
obtain the data $\{U_i, f_i\}$, where $U_i$ are an open cover of $Y$
and on each $U_i$ multiplication by $f_i$ is an isomorphism 
$\bigoplus \mathcal{L}(-D_\chi) \xrightarrow{\sim} \bigoplus
\mathcal{L}(-D'_\chi)$. One can readily check that such $\{U_i, f_i\}$
must define a Cartier divisor and that the corresponding Weil divisor
is the requisite divisor $N$.
\end{proof}

\begin{cor} \label{cor-equiv-classes-normalized-sets}
In each equivalence classes of \gnat-families there is precisely
one family whose reductor set is normalised.
\end{cor}

\subsection{Reductor Condition}

  We now investigate when is a set $\{D_\chi\}$ of $G$-divisors 
a reductor set. 
  
  This issue is the issue of $\bigoplus \mathcal{L}( - D_\chi)$
actually being $(\twalg) \otimes \mathcal{O}_Y$-module. By definition
it is a sub-$\mathcal{O}_Y$-module of $K(\mathbb{C}^n)$, but 
there is no a priori reason for it to also be closed under the
natural $\twalg$-action on $K(\mathbb{C}^n)$. If it is closed,
it can be checked that it trivially satisfies all the other
requirements in Proposition $\ref{prps-tfae}$, item  
$\ref{item-embeds-into-kcn}$, which makes it a  $\gnat$-family. 
We further observe that $\bigoplus \mathcal{L}( - D_\chi)$ is always 
closed under the action of $G$, so it all boils down to the closure 
under the action of $\regring$. 

  Recall, that we write $R_G$ for $R \cap \kgc_$, the
$G$-homogeneous regular polynomials, and $R_\chi$ for $R \cap
\kgchi_$, the $G$-homogeneous regular polynomials of weight $\chi \in
G^\vee$. 

\begin{prps}[Reductor Condition] \label{prps-reductor-condition}
Let $\{D_\chi\}_{\chi \in G^\vee}$ be a set with each $D_\chi$ a 
$\chi$-Weil divisor. Then it is a reductor set if and only if, 
for any $f \in \regring_G$, the divisor 
\begin{align} \label{eqn-reductor-condition} 
D_\chi + (f) - D_{\chi \rho(f)} \geq 0
\end{align}
i.e. it is effective. 
\end{prps}

\bf Remarks: \rm 
\begin{enumerate}
\item If we choose a $G$-eigenbasis of $\givrep$, then its dual basis, 
a set of basic monomials $x_1, \dots, x_n$, generates $\regring_G$
as a semi-group. As condition \eqref{eqn-reductor-condition} is
multiplicative on $f$, it is sufficient to check it only for $f$ being 
one of $x_i$. This leaves us with a finite number of inequalities to check. 

\item Numerically, if we write each $D_\chi$ as $\sum q_{\chi,P} P$,
inequalities $\eqref{eqn-reductor-condition}$ subdivide into
independent sets of inequalities 
\begin{align} \label{eqn-red-cond-num} q_{\chi,P} +
v_P(f) - q_{\chi \rho(f), P} \geq 0 \quad\quad \forall\; \chi \in G^\vee
\end{align} a set for each prime divisor $P$ on $Y$.
This shows that a $\gnat$-family can be specified independently at each 
prime divisor of $Y$: we can construct reductor sets $\{D_\chi\}$ by 
independently choosing for each prime divisor $P$ any of the sets 
of numbers $\{q_{\chi, P}\}_{\chi \in G^\vee}$ which satisfy 
\eqref{eqn-red-cond-num}.

\item There is an interesting link here with the work of Craw,
Maclagan and Thomas in \cite{Craw-Maclagan-Thomas-05-I} which 
bears further investigation. In a toric context, they have 
rediscovered these inequalities as dual, in a certain sense, 
to the defining equations of the coherent component $Y_\theta$ of 
the moduli space $M_\theta$ of $\theta$-semistable $G$-constellations.  
They then use them to compute the distinguished $\theta$-semistable 
$G$-constellations parametrised by torus orbits of $Y_\theta$. In
particular, their Theorem 7.2 allows them to explicitly
write down the tautological gnat-family on $Y_\theta$ and
suggests that, up to a reflection, it is the gnat family which 
minimizes $\theta . \{ D_\chi \}$. We shall see an example
of that for the case of $Y_\theta = \hilb^G$ in our Proposition
\ref{prps-maxshift-ghilb}.
\end{enumerate}

\begin{proof} Take an open cover $U_i$ on which all 
$\mathcal{L}(-D_\chi)$ are trivialised and write $g_{\chi,i}$ for
the generator of $\mathcal{L}(-D_\chi)$ on $U_i$. As $\regring$ is
a direct sum of its $G$-homogeneous parts, it is sufficient to check
the closure under the action of just the homogeneous functions. Thus
it suffices to establish that  
for each $f \in R_G$, each $U_i$ and each $\chi \in G^\vee$
\begin{align*} f g_{\chi,i} \in \mathcal{O}_Y(U_i)
g_{\chi \rho(f),i} \end{align*} 

	On the other hand, with the notation above, $G$-Cartier
divisor $D_\chi + (f) - D_{\chi \rho(f)}$ is given on $U_i$ by
$\frac{f g_{\chi,i} }{g_{\chi\rho(f),i}}$ and it being effective is
equivalent to \begin{align*} \frac{f g_{\chi,i} }{g_{\chi \rho(f),i}}
\in \mathcal{O}_Y(U_i) \end{align*} for all $U_i$'s. The result follows.
\end{proof}

\subsection{Canonical family} \label{chapter-canonical}

 We have not yet given any evidence of any \gnat-families actually
existing on an arbitrary resolution $Y$ of $X$. 

\begin{prps}[Canonical family] \label{prps-canonical-family}

    Let $Y$ be a resolution of $X = \mathbb{C}^n/G$. Define the set 
$\{ C_\chi \}_{\chi \in G^\vee}$ of $G$-Weil divisors
by $$C_\chi = \sum v(P, \chi) P$$ 
where $P$ runs over all prime Weil divisors on $Y$ and $v(P, \chi)$ 
are the numbers introduced in Definition 
$\ref{defn-homogeneous-valuations}$ (lifted to $[0,1) \subset
\mathbb{Q}$). 

Then $\{ C_\chi \}_{\chi \in G^\vee}$ is a reductor set.

We call the corresponding family 
\tt the canonical \gnat-family on $Y$\it.  
\end{prps} 

\begin{proof} 
We must show that $\{C_\chi\}$ satisfies the inequalities 
$\eqref{eqn-reductor-condition}$.
Choose any $\chi \in G^\vee$, any $f \in \regring_G$ 
and any prime divisor $P$ on $Y$.
Observe that $0 \leq v_P(\chi), v_P(\chi \rho(f)) <
1$ by definition,  while $v_P(f) \geq 0$ since $f^{|G|}$ is regular 
on all of $Y$. So we must have 
$$ v_P(\chi) + v_P(f) - v_P (\chi \rho(f)) > -1 $$ 
As the above expression must be integer-valued, we further have
\begin{align*} 
v_P(\chi) + v_P(f) - v_P(\chi \rho(f)) \geq 0 
\end{align*}
as required. \end{proof}

This family has a following geometrical description:
\begin{prps}
On any resolution $Y$, the canonical family is isomorphic to 
the pushdown to $Y$ of the structure sheaf $\mathcal{N}$ of 
the normalisation of the reduced fiber product 
$Y \times_X \mathbb{C}^n$.
\end{prps}
\begin{proof}
First we construct a $(\twalg) \otimes \mathcal{O}_Y$-module
embedding of $\mathcal{N}$ into $K(\mathbb{C}^n)$. Let $\alpha$ 
be the map $\mathcal{O}_Y \otimes_{\regring^G} \regring
\rightarrow K(\mathbb{C}^n)$ which sends $ a \otimes b$ to $ab$.
It is $\twalg \otimes \mathcal{O}_Y$-equivariant.
If we show that $\ker \alpha$ is the nilradical of 
$\mathcal{O}_Y \otimes_{\regring^G} \regring$, 
then $\mathcal{N}$ can be identified with the integral closure of 
the image of $\alpha$ in $K(\mathbb{C}^n)$. Due to 
$G$-equivariance $\alpha$ decomposes as $\bigoplus_{\chi \in
G^\vee} \alpha_\chi$ with each $\alpha_\chi$ a morphism
$\mathcal{O}_Y \otimes_{\regring^G} \regring_\chi \rightarrow
K_\chi(\mathbb{C}^n)$. Observe that 
$ (\mathcal{O}_Y \otimes_{\regring^G} \regring_\chi)^{|G|}
\subset \mathcal{O}_Y \otimes_{\regring^G} \regring_\trch_ = 
\mathcal{O}_Y $
as a product of $|G|$ homogeneous functions is invariant. Hence
$ (\ker \alpha_\chi)^{|G|} \subset \ker \alpha_\trch_ = 0 $ as required.

 Write $\bigoplus_{\chi \in G^\vee} \mathcal{N}_\chi$ for the decomposition of
$\mathcal{N}$ into $G$-eigensheaves. Fix a point $p \in Y$
and observe that $f \in K_\chi(\mathbb{C}^n)$ is integral over 
the local ring $\mathcal{N}_p$ if and only if 
$f^{|G|} \in (\mathcal{N}_\trch_)_p = \mathcal{O}_{Y,p}$. 
Therefore 
$$ (\mathcal{N}_\chi)_p = \{ f \in K_\chi(\mathbb{C}^n) \;|\; \text{$G$-Weil
divisor }\divr(f) \text{ is effective at } p \} $$
In particular, the generator $c_\chi$ of $C_\chi$ at 
$p$ lies in $(\mathcal{N}_\chi)_p$. Observe further that for 
any $f \in (\mathcal{N}_\chi)_p$ the Weil divisor $\divr(f) - C_\chi$
is effective at $p$ as the coefficients of $C_\chi$ are just the fractional 
parts of those of $\divr(f)$ and the latter is effective. Therefore 
$c_\chi$ generates $(\mathcal{N}_\chi)_p$ as
$\mathcal{O}_{Y,p}$-module, giving $\mathcal{N}_\chi = \mathcal{L}(-C_\chi)$ 
as required. 
\end{proof}

\subsection{Symmetries} \label{section-symmetries}

Having demonstrated that the set of equivalence classes of \gnat-families 
is always non-empty, we now establish two types of symmetries which this 
set possesses. It is worth noting that from the description of
the symmetries of the chambers in the parameter 
space of the stability conditions for $G$-constellations described 
in \cite{Craw-Ishii-02}, Section 2.5, it follows that all the symmetries
described below take the subset of $\gnat$-families on $Y$ consisting 
of universal families of stable $G$-constellations into itself. 

\begin{prps}[Character Shift]
Let $\{D_\chi\}$ be a normalised reductor set. Then for any $\chi$ in 
$G^\vee$ 
\begin{align}
D'_{\chi  \lambda} = D_{\chi} - D_{\lambda^{-1}}
\end{align}
is also a normalised reductor set. We call it \tt \;the $\chi$-shift \rm of 
$\{D_\chi\}$.
\end{prps}
\begin{proof}
Writing out the reductor condition $\eqref{eqn-reductor-condition}$
for the new divisor set $\{D'_\chi\}$ we get:
$$ (D_\chi - D_{\lambda^{-1}}) + (m) - (D_{\chi \rho(m)} -
D_{\lambda^{-1}}) \geq 0 $$
Cancelling out $D_\lambda^{-1}$, we obtain precisely the reductor
condition for the original set $\{D_\chi\}$. And since 
$$ D'_\trch_
= D'_{\lambda^{-1} \lambda} = D_{\lambda^{-1}} - D_{\lambda^{-1}} = 0
$$
we see that the new reductor set is normalised.
\end{proof}
\bf NB: \rm Observe, that for a reductor set $\{D_\chi\}$
and for any $\chi$-Weil divisor $N$, the set $\{D_\chi + N\}$ is 
linearly equivalent to the $\chi$-shift of $\{D_\chi\}$.

\begin{prps}[Reflection] 
Let $\{D_\chi\}$ be a normalised reductor set. 
Then the set $\{ - D_{\chi^{-1}}\}$ is also a normalised reductor set, 
which we call \tt \;the reflection \rm of $\{D_\chi\}$.
\end{prps}
\begin{proof}
We need to show that 
\begin{align*}
- D_{\chi^{-1}} + (m) - ( - D_{\chi^{-1}\rho(m)^{-1}}) \geq 0
\end{align*}
Rearranging we get
\begin{align*}
D_{\chi^{-1} \rho(m)^{-1}} + (m) - D_{\chi^{-1} \rho(m)^{-1} \rho(m)}
\geq 0
\end{align*}
which is one of the reductor equations the original set $\{D_\chi\}$
must satisfy. As $D'_\trch_ = - D_\trch_ = 0$, the new set is
normalised.
\end{proof}

\subsection{Maximal shift family and finiteness}

  We now examine the individual line bundles $\mathcal{L}(-D_\chi)$ in 
a $\gnat$-family and show that the reductor condition 
imposes a restriction on how far apart from each other they can be. 

\begin{lemma} \label{lemma-gen-max-shift}
 Let $\{D_\chi\}$ be a reductor set. Write each $D_\chi$ as 
$\sum q_{\chi, P} P$, where $P$ ranges over all the prime Weil
divisors on $Y$. For any $\chi_1, \chi_2 \in G^\vee$ and for any prime 
Weil divisor $P$, we necessarily have 
\begin{align}
\min_{f \in \regring_{\chi_1/\chi_2}} v_P(f)  
\quad \geq \quad q_{\chi_1,P} - q_{\chi_2,P}
\quad \geq \quad 
- \min_{f \in \regring_{{\chi_2}/{\chi_1}}} v_P(f)
\end{align}
\end{lemma}
\begin{proof}
Both inequalities follow directly from the reductor condition 
\eqref{eqn-reductor-condition}: the right inequality by setting 
$\chi = \chi_1 \in G^\vee$, $\rho(f) = \frac{\chi_2}{\chi_1}$ 
and letting $f$ vary within $\regring_{\rho(f)}$; the left inequality
by setting $\chi = \chi_2$ and $\rho(f) = \frac{\chi_1}{\chi_2}$.
\end{proof}

This suggests the following definition:
\begin{defn} \label{def-maxshift}
For each character $\chi \in G^\vee$, we define the 
\tt maximal shift $\chi$-divisor $M_\chi$ \it to be
\begin{align} \label{eqn-max-shift-defn}
M_\chi = \sum_P (\min_{f \in \regring_\chi} v_P(f)) P
\end{align}
where $P$ ranges over all prime Weil divisors on $Y$.  
\end{defn}

\begin{lemma} \label{lemma-max-shift-is-reductor}
The $G$-Weil divisor set $\{M_\chi\}$ is a normalised reductor set. 
We call the corresponding family \tt the maximal shift $\gnat$-family
\rm on $Y$.
\end{lemma}
\begin{proof}
We need to show that for any $f \in \regring_G$ 
and any prime divisor $P$ 
$$ v_P(m_\chi) + v_P(f) - v_P(m_{\chi \rho(f)}) \geq 0 $$
where $m_\chi$ and $m_{\chi \rho(f)}$ are chosen to achieve the 
minimality in \eqref{eqn-max-shift-defn}. 

 Observe that $m_\chi f$ is also a $G$-homogeneous element of $\regring$,
therefore by the minimality of $v_P (m_{\chi \rho(f)})$ we have
$$ v_P(m_\chi f)  \geq  v_P (m_{\chi \rho(f)}) $$
as required. 

To establish that $M_\trch_ = 0$ we observe that for any 
$G$-homogeneous $f \in \regring$ we have $v_P(f) \geq 0$
on any prime Weil divisor $P$ as $f^|G|$ is globally regular.  
Moreover for $f$ in $\regring_\trch_ = \regring^G$ this lower bound
is achieved by $f = 1$. 
\end{proof}

Observe that with Lemma \ref{lemma-max-shift-is-reductor} we have
established another \gnat-family which always
exists on any resolution $Y$. While sometimes it coincides 
with the canonical family, generally the two are distinct. 

\begin{prps}[Maximal Shifts] \label{prps-max-shifts}
Let $\{D_\chi\}$ be a normalised reductor set.
Then for any $\chi \in G^\vee$  
\begin{align} \label{eqn-family-bounded}
M_\chi \geq
D_\chi \geq
- M_{\chi^{-1}}
\end{align}
Moreover both the bounds are achieved.
\end{prps}
\begin{proof}
To establish that \eqref{eqn-family-bounded} holds set $\chi_2
= \trch_$ in Lemma \ref{lemma-gen-max-shift}. Lemma 
\ref{lemma-max-shift-is-reductor} shows that the bounds are achieved.
\end{proof}

\begin{prps} \label{prps-maxfam-coeffs-not-zero}
If the coefficient of a maximal shift divisor $M_\chi$ at a prime
divisor $P \subset Y$ is non-zero, then either $P$ is an
exceptional divisor or the image of $P$ in $X$ is a branch divisor
of $\mathbb{C}^n \xrightarrow{q} X$.  
\end{prps}
\begin{proof}
Let $P$ be a prime divisor on $X$ which is not a branch divisor
of $q$. Let $\chi \in G^\vee$. By the defining formula
$\eqref{eqn-max-shift-defn}$ it suffices to find $f \in R_\chi$
such that $v_P(f) = 0$. 

As $\regring$ is a PID, there exist $t_1, \dots,
t_k \in \regring$ such that $(t_1), \dots, (t_k)$ are all the distinct
prime divisors lying over $P$ in $\mathbb{C}^n$. Observe
that the product $t_1 \dots t_n$ must be $G$-homogeneous. 
Since $P$ is not a branch divisor, there exists $u \in \regring$ such
that $t_1 \dots t_k u$ is invariant and $u \notin (t_i)$ for all $i$.
Then $u' = u^{\left|G\right| - 1}$ is a $G$-homogeneous function of
the same weight as $t_1 \dots t_k$ and $v_P(u') = 0$. Now take any $f
\in \regring_\chi$ and consider its factorization into irreducibles.
$G$-homogeneity of $f$ implies that all $t_i$ occur with the same
power $k$. Now replacing $(t_1 \dots t_n)^k$ in the factorization by
$(u')^k$ we obtain an element of $\regring_\chi$ whose valuation at 
$P$ is zero. 
\end{proof}

\begin{cor} \label{cor-equiv-classes-finite}
The number of equivalence classes of $\gnat$-families on $Y$ is
finite.
\end{cor}
\begin{proof}
Let $\{D_\chi\}$ be a normalised reductor set. Coefficients 
of $D_\chi$ at prime divisors $P$ of $Y$ have fixed fractional
parts (Definition \ref{defn-gweil-divisor}), are bound above and 
below (Proposition \ref{prps-max-shifts}) and are zero at all but 
finite number of $P$ (Proposition \ref{prps-maxfam-coeffs-not-zero}). 
This leaves only a finite number of possibilities.   
\end{proof}

For one particular resolution $Y$ the family provided by 
the maximal shift divisors has a nice geometrical description. 

\begin{prps} \label{prps-maxshift-ghilb}
Let $Y = \hilb^G \mathbb{C}^n$, the coherent component of 
the moduli space of $G$-clusters in $\mathbb{C}^n$. If $Y$ is smooth, 
then $\bigoplus \mathcal{L}(-M_\chi)$ is the universal 
family $\mathcal{F}$ of $G$-clusters parametrised by
$Y$, up to the usual equivalence of families. 
\end{prps}

\begin{proof}
Firstly $\mathcal{F}$ is a \gnat-family, as
over any set $U \subset X$ such that $G$ acts freely on $q^{-1}(U)$
we have $\mathcal{F} |_U \simeq \pi^* q_* 
\mathcal{O}_{\mathbb{C}^n} |_U$.
Write $\mathcal{F}$ as $\oplus \mathcal{L}(-D_\chi)$ for some
reductor set $\{ D_\chi \}$. Take an open cover $\{U_i\}$ of $Y$ and
consider the generators $\{f_{\chi,i}\}$ of $D_\chi$ on each
$U_i$. Working up to equivalence, we can consider $\{D_\chi\}$ 
to be normalised and so $f_{\trch_, i} = 1$ for all $U_i$.

Now any $G$-cluster $Z$ is given by some invariant ideal $I \subset
\regring$ and so the corresponding $G$-constellation
$H^0(\mathcal{O}_Z)$
is given by $\regring / I$. In particular note that $R/I$ is generated
by $R$-action on the generator of $\trch_$-eigenspace. Therefore
any $f_{\chi, i}$ is generated from $f_{\trch_, i} = 1$ by
$\regring$-action, which means that all $f_{\chi, i}$ lie in
$\regring$.

But this means that for any prime Weil divisor $P$ on $Y$ we have
\begin{align*}
v_P(f_{\chi, i}) \geq \min_{f \in \regring_\chi} v_P(f)
\end{align*}
and therefore $D_\chi \geq M_\chi$. Now Proposition \ref{prps-max-shifts}
forces the equality. 
\end{proof}

\section{Conclusion}

We summarise the results achieved in the following theorem:

\begin{theorem}[Classification of $\gnat$-families]
\label{theorem-classification}

Let $G$ be a finite abelian subgroup of $\gl_n(\mathbb{C})$, $X$ the
quotient of $\mathbb{C}^n$ by the action of $G$, $Y$ nonsingular
and $\pi:\; Y \rightarrow X$ a proper birational map. 
Then isomorphism classes of \gnat-families on $Y$
are in $1$-to-$1$ correspondence with linear equivalence classes 
of $G$-divisor sets $\{D_\chi\}_{\chi \in G^\vee}$, each 
$D_\chi$ a $\chi$-Weil divisor, which satisfy the inequalities
\begin{align*}
D_\chi + (f) - D_{\chi \rho(f)} \geq 0 \quad \forall\; \chi \in
G^\vee, \text{$G$-homogeneous } f \in \regring
\end{align*}
Such a divisor set $\{D_\chi\}$ corresponds then to a 
$\gnat$-family $\bigoplus \mathcal{L}(-D_\chi)$.

This correspondence descends to a $1$-to-$1$ correspondence between 
equivalence classes of \gnat-families and sets $\{D_\chi\}$ 
as above and with $D_{\chi_0} = 0$. Furthermore, each divisor $D_\chi$
in such a set satisfies inequality 
\begin{align*} M_\chi \geq
D_\chi \geq
- M_{\chi^{-1}}
\end{align*}
where $\{ M_{\chi} \}$ is a fixed divisor set defined by
\begin{align*}
M_\chi = \sum_P (\min_{f \in \regring_\chi} v_P(f)) P
\end{align*}
As a consequence, the number of equivalence classes of $\gnat$-families 
is finite. 
\end{theorem}

\begin{proof}
Corollary $\ref{cor-iso-families-reductors-lineq}$
establishes the correspondence between isomorphism classes of 
$\gnat$-families and linear equivalence classes of reductor sets.
Proposition $\ref{prps-reductor-condition}$ gives description 
of reductor sets as the divisor sets satisfying the reductor
condition inequalities.

Corollary \ref{cor-equiv-classes-normalized-sets} gives the
correspondence on the level of equivalence classes of $\gnat$-families
and normalised reductor sets. Proposition \ref{prps-max-shifts}
establishes the bounds on the set of all normalised reductor sets and 
Corollary \ref{cor-equiv-classes-finite} uses it to show that the
set of all normalised reductor sets is finite.
\end{proof}

\bibliography{../../references}

\providecommand{\bysame}{\leavevmode\hbox to3em{\hrulefill}\thinspace}
\providecommand{\MR}{\relax\ifhmode\unskip\space\fi MR }
\providecommand{\MRhref}[2]{%
  \href{http://www.ams.org/mathscinet-getitem?mr=#1}{#2}
}
\providecommand{\href}[2]{#2}
\begin{thebibliography}{ACvdE05}

\bibitem[ACvdE05]{Adjamagbo-Charbonnel-VanDenEssen-05}
Kossivi Adjamagbo, Jean-Yves Charbonnel, and Arno van~den Essen, \emph{On ring
  homomorphisms of {A}zumaya algebras}, preprint math.RA/0509188, (2005).

\bibitem[AG60]{Auslander-Goldamn-1960}
M.~Auslander and O.~Goldman, \emph{The {B}rauer group of a commutative ring},
  Trans. Amer. Math. Soc. \textbf{97} (1960), 367--409.

\bibitem[BKR01]{BKR01}
T.~Bridgeland, A.~King, and M.~Reid, \emph{The {M}c{K}ay correspondence as an
  equivalence of derived categories}, J. Amer. Math. Soc. \textbf{14} (2001),
  535--554, math.AG/9908027.

\bibitem[BO95]{BonOr95}
A.~Bondal and D.~Orlov, \emph{Semi-orthogonal decompositions for algebraic
  varieties}, preprint math.AG/950612, (1995).

\bibitem[Bou98]{Bouvier-1998}
Catherine Bouvier, \emph{Diviseurs essentiels, composantes essentielles des
  vari\'et\'es toriques singuli\'eres}, Duke Math J. \textbf{91 no.3} (1998),
  609--620.

\bibitem[Bri99]{Bridg97}
T.~Bridgeland, \emph{Equivalence of triangulated categories and
  {F}ourier-{M}ukai transforms}, Bull. London Math. Soc. \textbf{31} (1999),
  25--34, math.AG/9809114.

\bibitem[CI04]{Craw-Ishii-02}
A.~Craw and A.~Ishii, \emph{Flops of ${G}$-$\hilb$ and equivalences of derived
  category by variation of {GIT} quotient}, Duke Math J. \textbf{124} (2004),
  no.~2, 259--307, arXiv:math.AG/0211360.

\bibitem[CMT07a]{Craw-Maclagan-Thomas-05-I}
A.~Craw, D.~Maclagan, and R.R. Thomas, \emph{Moduli of {M}c{K}ay quiver
  representations {I}: the coherent component}, Proc. London Math. Soc.
  \textbf{95} (2007), no.~1, 179 -- 198.

\bibitem[CMT07b]{Craw-Maclagan-Thomas-05-II}
\bysame, \emph{Moduli of {M}c{K}ay quiver representations {II}: Grobner basis
  techniques}, J. Algebra \textbf{316} (2007), no.~2, 514--535,
  arXiv:math.AG/0611840.

\bibitem[Cra01]{Craw-thesis}
A.~Craw, \emph{The {M}c{K}ay correspondence and representations of the
  {M}c{K}ay quiver}, Ph.D. thesis, University of {W}arwick, 2001.

\bibitem[dB02]{vdBergh2002}
Michel~Van den Bergh, \emph{Non-commutative crepant resolutions}, The Legacy of
  Niels Hendrik Abel, Springer, 2002, pp.~749--770.

\bibitem[Gar86]{Garl86}
D.J.H. Garling, \emph{A course in {G}alois theory}, Cambridge University Press,
  1986.

\bibitem[Har77]{Harts77}
R.~Hartshorne, \emph{Algebraic geometry}, Springer-Verlag, 1977.

\bibitem[Kal08]{Kaledin-05}
D.~Kaledin, \emph{Derived equivalences by quantization}, Geom. Funct. Anal.
  \textbf{17} (2008), no.~6, 1968--2004, arXiv:math.AG/0504584.

\bibitem[Kaw05]{Kawamata-LogCrepantBirationalMapsAndDerivedCategories}
Y.~Kawamata, \emph{Log crepant birational maps and derived categories}, J.
  Math. Sci. Univ. Tokyo \textbf{12} (2005), no.~2, 211--231, math.AG/0311139.

\bibitem[KI03]{Kollar-Ishii-2003}
J\'anos Koll\'ar and Shihoko Ishii, \emph{The {N}ash problem on arc families of
  singularities}, Duke Math J. \textbf{120 no.3} (2003), 601--620.

\bibitem[Kin94]{King94}
A.~King, \emph{Moduli of representations of finite-dimensional algebras},
  Quart. J. Math. Oxford \textbf{45} (1994), 515--530.

\bibitem[Log04]{Logvinenko-thesis}
T.~Logvinenko, \emph{Families of ${G}$-{C}onstellations parametrised by
  resolutions of quotient singularities}, Ph.D. thesis, University of {B}ath,
  2004.

\bibitem[Log08]{Logvinenko-DerivedMcKayCorrespondenceViaPureSheafTransforms}
\bysame, \emph{Derived {M}c{K}ay correspondence via pure-sheaf transforms},
  Math. Ann. \textbf{341} (2008), no.~1, 137--167.

\bibitem[Mil80]{Milne-EtaleCohomology}
J.S. Milne, \emph{{\'E}tale cohomology}, Princeton University Press, 1980.

\bibitem[Rei97]{Kinosaki-97}
M.~Reid, \emph{{M}c{K}ay correspondence}, preprint math.AG/9702016, (1997).

\bibitem[Ver00]{Verbitsky-2000}
Misha Verbitsky, \emph{Holomorphic symplectic geometry and orbifold
  singularities}, Asian J. Math. \textbf{4} (2000), 553--563.

\bibitem[Zar58]{Zariski-1958}
Oscar Zariski, \emph{On the purity of the branch locus of algebraic functions},
  Proc. N. A. S. \textbf{44} (1958), 791--796.

\end{thebibliography}
\bibliographystyle{amsalpha}

\begin{tabular}{l l}
E-mail: & \tt T.Logvinenko@liv.ac.uk \rm \\
Address: & \tt Department of Mathematical Sciences, \rm \\
&\tt University of Liverpool, \rm\\
&\tt Peach Street, \rm\\
&\tt Liverpool, L69 7ZL, \rm\\
&\tt UK \rm 
\end{tabular}

\end{document}